\documentclass[12pt]{amsart}
\numberwithin{equation}{section}
\newtheorem{theorem}{Theorem}
\numberwithin{theorem}{section}

\numberwithin{theorem}{section} \numberwithin{lemma}{section}

\numberwithin{definition}{section}

\numberwithin{corollary}{section}

\numberwithin{remark}{section}

\numberwithin{proposition}{section}
\def\b{\begin{equation}}
\def\e{\end{equation}}
\newcommand{\ignore}[1]{}
\normalsize
\date {April 09, 2007}
\begin{document}
\pagenumbering{arabic} \pagenumbering{arabic}\setcounter{page}{1}
\tracingpages 1
\title{Hardy and Rellich type inequalities with remainders for Baouendi-Grushin vector fields}
\author{Ismail Kombe}
\dedicatory {}
\address{Ismail Kombe, Mathematics Department\\ Dawson-Loeffler Science
\&Mathematics Bldg\\
Oklahoma City University \\
2501 N. Blackwelder, Oklahoma City, OK 73106-1493}
\email{ikombe@okcu.edu}
\begin{abstract}
In this paper we study Hardy and Rellich type inequalities  for
Baouendi-Grushin vector fields : $\nabla_{\gamma}=(\nabla_x,
|x|^{2\gamma}\nabla_y)$ where  $\gamma>0$, $\nabla_x$ and
$\nabla_y$ are usual gradient operators in the variables $x\in
\mathbb{R}^m$ and $y\in\mathbb{R}^k$, respectively. In the first
part of the paper, we prove some weighted Hardy type inequalities
with remainder terms. In the second part, we prove two versions of
weighted Rellich type inequality on the whole space. We find sharp
constants for these inequalities. We also obtain their improved
versions for bounded domains.
\end{abstract}
\maketitle
\section{Introduction}
This paper is concerned with Hardy and Rellich type inequalities
with remainder terms for Baouendi-Grushin vector fields. Let $x\in
\mathbb{R}^m$, $y\in \mathbb{R}^k$, $\gamma>0$ and $n=m+k$, with
$m, k\ge 1$. Then the following Hardy type inequality for
Baouendi-Grushin vector fields has been proved by Garofalo
\cite{Garofalo},
\begin{equation}
\int_{
\mathbb{R}^n}\big(|\nabla_{x}\phi|^2+|x|^{2\gamma}|\nabla_{y}\phi|^2\big)dxdy
 \ge \Big(\frac{Q-2}{2}\Big)^2\int_{\mathbb{R}^n} (\frac{|x|^{2\gamma}}{|x|^{2+2\gamma}+(1+\gamma^2)^2|y|^2})\phi^2
 dxdy
\end{equation}
where $\phi\in C_0^{\infty}(\mathbb{R}^m\times
\mathbb{R}^k\setminus \{(0,0)\})$ and $Q=m+(1+\gamma)k$. Here,
$\nabla_x\phi$ and $\nabla_y\phi$ denotes the gradients of $\phi$
in the variables $x$ and $y$, respectively. A similar inequality
with the same sharp constant $(\frac{Q-2}{2})^2$ holds if $
\mathbb{R}^n$ replaced by $\Omega$ and $\Omega$ contains the
origin \cite{D'Ambrosio}. If $\gamma=0$ then it is clear that the
inequality (1.1) recovers the classical Hardy inequality in
$\mathbb{R}^n$
\begin{equation}
\int_{\mathbb{R}^n}|\nabla\phi(z)|^2dz\ge
\Big(\frac{n-2}{2}\Big)^2\int_{\mathbb{R}^n}
\frac{|\phi(z)|^2}{|z|^2}dz
\end{equation}
where $z=(x,y)\in\mathbb{R}^m\times \mathbb{R}^k$ and the constant
$(\frac{n-2}{2})^2$ is sharp.  There exists a large literature
concerning with the Hardy inequalities and, in particular, sharp
inequalities as well as their improved versions which have
attracted a lot of attention because of their application to
singular problems (See \cite{Baras-Goldstein},
\cite{Peral-Vazquez}, \cite{Brezis-Vazquez}, \cite{Garcia-Peral},
\cite{Cabre-Martel}, \cite{Vazquez-Zuazua}, \cite{Kombe1} and
references therein).

A sharp improvement of the Hardy inequality (1.2) was discovered
by Brezis and V\'azquez \cite{Brezis-Vazquez}. They proved that
for a bounded domain $\Omega\subset \mathbb{R}^n$
 \begin{equation}
\int _{\Omega} |\nabla \phi(z)|^2dz\ge
\Big(\frac{n-2}{2}\Big)^2\int_{\Omega}
\frac{|\phi(z)|^2}{|z|^2}dz+\mu\big(\frac{\omega_n}{|\Omega|}\big)^{2/n}\int
_{\Omega} \phi^2dz,
 \end{equation}
where $\phi\in C_0^{\infty} (\Omega)$, $\omega_n$  and $|\Omega|$
denote the $n$-dimensional Lebesgue measure of the unit ball
$B\subset \mathbb{R}^n$  and the domain $\Omega$ respectively.
Here $\mu$ is the first eigenvalue of the Laplace operator in the
two dimensional unit disk and it is optimal when $\Omega$ is a
ball centered at the origin. In a recent paper Abdelloui, Colorado
and Peral \cite{ACP} obtained, among other things, the following
improved Caffarelli-Kohn-Nirenberg inequality
\begin{equation}
\int _{\Omega} |\nabla \phi(z)|^2|z|^{-2a}dz\ge
\Big(\frac{n-2a-2}{2}\Big)^2\int_{\Omega}
\frac{|\phi(z)|^2}{|z|^{2a+2}}dz+C\Big(\int _{\Omega}
|\nabla\phi|^q|z|^{-aq}\Big)^{2/q}dz
\end{equation}
where $\phi\in C_0^{\infty} (\Omega)$, $-\infty<a<\frac{n-2}{2}$,
$1<q<2$ and $C=C(q,n,\Omega)>0$. Motivated by these results, our
first goal is to find improved weighted Hardy type inequalities
for Baouendi-Grushin vector fields.

It is well known that an important extension of Hardy's inequality
to higher-order derivatives is the following Rellich inequality
\begin{equation} \int_{\mathbb{R}^n}|\Delta \phi(z)|^2dz\ge
\frac{n^2(n-4)^2}{16}\int_{\mathbb{R}^n
}\frac{|\phi(z)|^2}{|z|^4}dz\end{equation} where $\phi\in
C_0^{\infty}(\mathbb{R}^n\setminus \{0\}$), $n \neq 2$ and the
constant $\frac{n^2(n-4)^2}{16}$   is sharp.  Davies and Hinz
\cite{Davies-Hinz}, among other results, obtained sharp weighted
Rellich inequalities of the form
\begin{equation} \int_{\mathbb{R}^n}\frac{|\Delta\phi(z)|^2}{|z|^{\alpha}}
dz\ge C\int_{
\mathbb{R}^n}\frac{|\phi(z)|^2}{|z|^{\beta}}dz\end{equation}  for
suitable values of $\alpha, \beta, p$ and $\phi\in
C_0^{\infty}(\mathbb{R}^n\setminus \{0\}$). In a recent paper,
Tertikas and Zographopoulos \cite{TZ}, among other results,
obtained the following new Rellich type inequalities that connects
first to second order derivatives:
\begin{equation}
\int_{\mathbb{R}^n}|\Delta \phi|^2dz\ge
\frac{n^2}{4}\int_{\mathbb{R}^n}\frac{|\nabla \phi|^2}{|z|^2}dz
\end{equation}
where $\phi\in C_0^{\infty}(\mathbb{R}^n\setminus \{0\})$ and  the
constant $\frac{n^2}{4}$ is sharp.  Recently, Kombe \cite{Kombe2}
obtained analogues of (1.6) and (1.7), and their improved versions
on Carnot groups. Motivated by the above results, our second goal
is to find sharp weighted Rellich type inequalities and their
improved versions for Baouendi-Grushin vector fields in that they
do not arise from any Carnot group. We should also mention that
Kombe and \"Ozaydin \cite{Kombe-Ozaydin} obtained (under some
geometric assumptions) improved Hardy and Rellich inequalities on
a Riemannian manifold that does not recover our current results.
Analogue inequalities for the Greiner vector fields will be given
in a forthcoming paper \cite{Kombe3}.

\section{Notations and Back ground material}

In this section, we shall collect some notations, definitions and
preliminary facts which will be used throughout the article. The
generic point is $z=(x_1,..., x_m, y_1,...,y_k)=(x,y)\in
\mathbb{R}^m\times \mathbb{R}^k$  with $m,k\ge 1, m+k=n$. The
sub-elliptic gradient is the $n$ dimensional vector field given by
\begin{equation}
\nabla_{\gamma}=(X_1, \cdots, X_m, Y_1, \cdots, Y_k)
\end{equation}
where
\begin{equation}
X_j=\frac{\partial}{\partial x_j}, \quad j=1, \cdots, m , \quad
Y_j=|x|^{\gamma}\frac{\partial}{\partial y_j}, \quad j=1, \cdots,
k.\end{equation} The Baouendi-Grushin operator  on
$\mathbb{R}^{m+k}$  is the operator
\begin{equation}\Delta_{\gamma}=\nabla_{\gamma}\cdot\nabla_{\gamma}
=\Delta_x+|x|^{2\gamma}\Delta_y,
\end{equation} where $\Delta_x$ and $\Delta_y$ are Laplace
operators in the variables $x\in \mathbb{R}^m$ and $y\in
\mathbb{R}^k$, respectively (see \cite{Baouendi}, \cite{Grushin1},
\cite{Grushin2}). If $\gamma$ is an even positive integer then
$\Delta_{\gamma}$ is a sum of squares  of $C^{\infty}$ vector
fields satisfying H\"ormander finite rank condition: \textit{rank
Lie} $ [\ X_1, \cdots, X_m, Y_1, \cdots, Y_k ]\ =n.$ The
anisotropic dilation attached to $\Delta_{\gamma}$ is given by
\[\delta_{\lambda}(z)=(\lambda x, {\lambda}^{\gamma+1}y), \quad
\lambda>0, \quad z=(x,y)\in \mathbb{R}^{m+k}.\] The change of
variable formula for the Lebesgue measure gives that\[d\circ
\delta_{\lambda}(x,y)=\lambda^{Q}dxdy,\] where
\[Q=m+(1+\gamma)k\] is the homogeneous dimension with respect to
dilation $\delta_{\lambda}$. For $z=(x,y)\in \mathbb{R}^m\times
\mathbb{R}^k$, let
\begin{equation}
\rho=\rho(z):=\Big(|x|^{2(1+\gamma)}+(1+\gamma)^2|y|^2\Big)^{\frac{1}{2(1+\gamma)}}.\end{equation}
By direct computation we get
\[
|\nabla_{\gamma}\rho|=\frac{|x|^{\gamma}}{\rho^{\gamma}}.
\]
Let $f\in C^2(0, \infty)$ and define  $u=f(\rho)$ then we have the
following useful formula
\begin{equation}
\Delta_{\gamma}u=\frac{|x|^{2\gamma}}{\rho^{2\gamma}}\Big(f''+\frac{Q-1}{\rho}f'\Big).
\end{equation}
We let $B_{\rho}=\{z\in \mathbb{R}^{n}\mid \rho(z)<r\}$,
$B_{\tilde\rho}=\{z\in \mathbb{R}^{n} \mid \tilde\rho(z,0)<r\}$
and call these sets, respectively,  $\rho$-ball and
Carnot-Carath\'eodory metric ball centered at the origin with
radius $r$. The Carnot-Carath\'eodory distance $\tilde\rho$
between the points $z$ and$z_0$ is defined by
\[\tilde\rho(z,z_0)=\text{inf}\{\text{length}(\eta)\,|\,\eta\in \mathcal{K}\}\] where the set $\mathcal{K}$ is the set of all curves
$\eta$ such that $\eta(0)=z$, $\eta(1)=z_0$ and $\dot{\eta}(t)$ is
in \textit{span}$\{X_1(\eta(t)),..., X_m(\eta(t)), Y_1(\eta(t)),
..., Y_k(\eta(t)) \}$. If $\gamma$ is a positive even integer then
Carnot-Carath\'eodory distance of $z$ from the origin $\tilde
\rho(z,0)$ is comparable to $\rho(z)$.  ( See \cite{Franchi} and
\cite{Beliche} for further details.)
\medskip

It is well known that  Sobolev and Poincar\'e type inequalities
are important in the study of partial differential equations,
especially in the study of those arising from geometry and
physics. In\cite{Franchi}, Franchi, Gutierrez and Wheeden obtained
the following  Sobolev-Poincar\'e inequality for metric balls
associated with Baouendi-Grushin type operators:
\begin{equation}
\Big(\frac{1}{w_1(B)}\int_{B}|\nabla_{\gamma}\phi|^p
w_1(z)dz\Big)^{1/p}\ge \frac{1}{cr}
\Big(\frac{1}{w_2(B)}\int_{B}|\phi(z)|^q w_2(z)dz\Big)^{1/q}
\end{equation} where $\phi\in C_0^{\infty}(B)$ and the weight functions $w_1$ and $w_2$  satisfies
some certain conditions. Here, $c$ is independent of $\phi$ and
$B$, $1\le p\le q<\infty$ and $w(B)=\int_B w(z)dz$. If $w_1=w_2=1$
then Monti \cite{Monti} obtained the  following sharp Sobolev
inequality
\begin{equation}
\Big(\int_{\mathbb{R}^n}\big(|\nabla_x \phi|^2+
|x|^{2\gamma}|\nabla_y \phi|^2\big)dxdy\Big)^{1/2}\ge C
\Big(\int_{\mathbb{R}^n}|\phi|^{\frac{2Q}{Q-2}}dxdy\Big)^{\frac{Q-2
}{2Q}}
\end{equation} where $C=C(m,k,\alpha)>0$.
\medskip

\section{Improved Hardy-type inequalities}
In this section we study improved Hardy type inequalities. These
inequalities plays key role in establishing improved Rellich type
inequalities. In the various integral inequalities below (Section
3 and Section 4), we allow the values of the integrals on the
left-hand sides to be $+\infty$. The following theorem is the
first result of this section.

\begin{theorem} Let $\gamma$ be an even positive integer, $\alpha \in \mathbb{R}$,  $-\frac{m}{\gamma}<t<\frac{m}{\gamma}$,
 and $Q+\alpha-2>0$. Then the following inequality is valid
\begin{equation}\begin{aligned}\int_{B_{\rho}} \rho^{\alpha}|\nabla_{\gamma} \rho|^t|\nabla_{\gamma}
\phi|^2dz &\ge \Big(\frac{Q+\alpha-2}{2}\Big)^2 \int_{B_{\rho}}
\rho^{\alpha}\frac{|\nabla_{\gamma} \rho|^{t+2}}{\rho^2}\phi
^2dz\\& +\frac{1}{C^2r^2}\int _{B_{\rho}}
\rho^{\alpha}|\nabla_{\gamma} \rho|^t\phi^2dz
\end{aligned}\end{equation}
for all compactly supported smooth function $\phi\in
C_0^{\infty}(B_{\rho})$.
\end{theorem}
\proof

 Let $\phi=\rho^{\beta}\psi\in
C_0^{\infty}(B_{\rho})$ and $\beta\in \mathbb{R}\setminus \{0\}$.
A direct calculation
 shows that
 \begin{equation}
 \begin{aligned}
\int _{B_{\rho}} \rho^{\alpha}|\nabla_{\gamma}
\rho|^t|\nabla_{\gamma} \phi|^2dz &= \beta^2 \int
_{B_{\rho}}\rho^{\alpha+2\beta-2}|\nabla_{\gamma} \rho|^{t+2}\psi^2 dz\\
&+ 2\beta\int _{B_{\rho}} \rho^{\alpha+2\beta-1}|\nabla_{\gamma}
\rho|^t\psi\nabla_{\gamma} \rho\cdot \nabla_{\gamma} \psi dz
\\ &+ \int _{B_{\rho}}\rho^{\alpha+2\beta}|\nabla_{\gamma}
\rho|^t|\nabla_{\gamma}\psi|^2dz. \end{aligned}
\end{equation}
Applying integration by parts to the middle term and using the
following fact \[\nabla_{\gamma}\cdot\big(
\rho^{\alpha+2\beta-1}|\nabla_{\gamma} \rho|^t\nabla_{\gamma}
\rho\big)=
(Q+\alpha+2\beta-2)\rho^{\alpha+2\beta-2}|\nabla_{\gamma}
\rho|^{t+2}\] yields
\begin{equation}
 \begin{aligned}
\int _{B_{\rho}} \rho^{\alpha}|\nabla_{\gamma}
\rho|^t|\nabla_{\gamma} \phi|^2dz = f(\beta)\int
_{B_{\rho}}\rho^{\alpha+2\beta-2}|\nabla_{\gamma}
\rho|^{t+2}\psi^2 dz + \int
_{B_{\rho}}\rho^{\alpha+2\beta}|\nabla_{\gamma}
\rho|^t|\nabla_{\gamma}\psi|^2dz\end{aligned}
\end{equation}
where $f(\beta)= -\beta^2-\beta(\alpha+Q-2)$. Note that $f(\beta)$
attains the maximum for $\beta=\frac{2-\alpha-Q}{2}$, and this
maximum is equal to $C_{H}=(\frac{Q+\alpha-2}{2})^2$. Therefore we
have the following
\begin{equation}
 \begin{aligned}
\int _{B_{\rho}} \rho^{\alpha}|\nabla_{\gamma}
\rho|^t|\nabla_{\gamma} \phi|^2dz = C_{H}\int
_{B_{\rho}}\rho^{\alpha-2}|\nabla_{\gamma} \rho|^{t+2}\phi^2 dz +
\int _{B_{\rho}}\rho^{2-Q}|\nabla_{\gamma}
\rho|^t|\nabla_{\gamma}\psi|^2dz.\end{aligned}
\end{equation}
It is easy to show that the
 weight functions $w_1=w_2=\rho^{2-Q}|\nabla_{\gamma}\rho|^t$
satisfies the Muckenhoupt $A_2$ condition for
$-\frac{m}{\gamma}<t<\frac{m}{\gamma}$. Therefore  weighted
Poincar\'e inequality holds (see \cite{Franchi}, \cite{Lu},
\cite{FGaW}) and we have

\[\begin{aligned}\int _{B_{\rho}}\rho^{2-Q}|\nabla_{\gamma}
\rho|^t|\nabla_{\gamma}\psi|^2dz &\ge \frac{1}{C^2r^2}\int
_{B_{\rho}}
\rho^{2-Q}|\nabla_{\gamma} \rho|^t\psi^2dz\\
&=  \frac{1}{C^2r^2}\int _{B_{\rho}} \rho^{\alpha}|\nabla_{\gamma}
\rho|^t\phi^2dz\end{aligned}\] where $C$ is a positive constant
and $r^2$ is the radius of the ball $B_{\rho}$.

We now obtain the desired inequality
\begin{equation}\begin{aligned}
\int _{B_{\rho}} \rho^{\alpha}|\nabla_{\gamma}
\rho|^t|\nabla_{\gamma} \phi|^2dz \ge C_{H}\int
_{B_{\rho}}\rho^{\alpha-2}|\nabla_{\gamma} \rho|^{t+2}\phi^2 dz
+\frac{1}{C^2r^2}\int _{B_{\rho}} \rho^{\alpha}|\nabla_{\gamma}
\rho|^t\phi^2dz.\end{aligned}
\end{equation}
\endproof

 Using the same method, we have the following weighted
Hardy inequality which has a logarithmic remainder term. Similar
results in the Euclidean setting can be found in \cite{FT},
\cite{Adimurthi}, \cite{Wang-Willem}, \cite{ACP}.
\begin{theorem} Let $\alpha \in \mathbb{R}$,  $t\in \mathbb{R}$, $Q+\alpha-2>0$. Then the
following inequality is valid
\begin{equation}\begin{aligned}\int_{B_{\rho}} \rho^{\alpha}|\nabla_{\gamma} \rho|^t|\nabla_{\gamma}
\phi|^2dz \ge C_{H} \int_{B_{\rho}}
\rho^{\alpha-2}|\nabla_{\gamma} \rho|^{t+2}\phi ^2dz
+\frac{1}{4}\int _{B_{\rho}} \rho^{\alpha-2}|\nabla_{\gamma}
\rho|^{t+2}\frac{\phi^2}{(\ln\frac{r}{\rho})^2}dz
\end{aligned}\end{equation}
for all compactly supported smooth function $\phi\in
C_0^{\infty}(B_{\rho})$. \proof We have the following result from
(3.4):
\begin{equation}
 \begin{aligned}
\int _{B_{\rho}} \rho^{\alpha}|\nabla_{\gamma}
\rho|^t|\nabla_{\gamma} \phi|^2dz = C_{H}\int
_{B_{\rho}}\rho^{\alpha-2}|\nabla_{\gamma} \rho|^{t+2}\phi^2 dz +
\int _{B_{\rho}}\rho^{2-Q}|\nabla_{\gamma}
\rho|^t|\nabla_{\gamma}\psi|^2dz.\end{aligned}
\end{equation}
Let $\varphi \in C_0^{\infty}(B_{\rho})$ and set $\psi(z)=(\ln
\frac{r}{\rho})^{1/2}\varphi(z)$. A direct computation shows that
\begin{equation}\begin{aligned}
\int _{B_{\rho}}\rho^{2-Q}|\nabla_{\gamma}
\rho|^t|\nabla_{\gamma}\psi|^2dz&\ge
\frac{1}{4}\int_{B_{\rho}}\rho^{-Q}|\nabla_{\gamma}
\rho|^{t+2}\frac{\psi^2}{(\ln\frac{r}{\rho})^2}dz
\\
&=\frac{1}{4}\int_{B_{\rho}}\rho^{\alpha-2}|\nabla_{\gamma}
\rho|^{t+2}\frac{\phi^2}{(\ln\frac{r}{\rho})^2}dz.
\end{aligned}\end{equation} Substituting (3.8) into (3.7) which yields
the desired inequality (3.6). \qed
\end{theorem}
\medskip

We now first prove the following weighted $L^p$-Hardy inequality
which plays an important role in the proof of Theorem 3.3, Theorem
4.1 and Theorem 4.5.
\begin{theorem} Let $\Omega$ be either bounded or unbounded domain with smooth
boundary which contains origin, or $\mathbb{R}^n$.  Let $\alpha
\in \mathbb{R}$, $t\in \mathbb{R}$, $ 1\le p<\infty$ and
$Q+\alpha-p>0$. Then the following inequality holds
\begin{equation}\int _{\Omega}\rho ^{\alpha}|\nabla_{\gamma}\rho|^t|\nabla_{\gamma}\phi|^p dz\ge
\Big(\frac{Q+\alpha-p}{p}\Big)^p\int _{\Omega}
\rho^{\alpha}|\nabla_{\gamma}\rho|^t\frac{|\nabla_{\gamma}\rho|^p}{\rho^
p}|\phi|^pdz\end{equation} for all compactly supported smooth
functions $\phi\in C_0^{\infty}(\Omega)$.

\proof Let $\phi=\rho^{\beta}\psi\in C_0^{\infty}(\Omega)$ and
$\beta\in \mathbb{R}-\{0\}$. We have
\[
|\nabla_{\gamma}(\rho^{\beta}\psi)|=|\beta
\rho^{\beta-1}\psi\nabla_{\gamma} \rho+
\rho^{\beta}\nabla_{\gamma} \psi|.
\]
We now use the following inequality which is  valid for any $a, b
\in \mathbb{R}^n$ and  $p>2$,

\[
|a+b|^p-|a|^p\ge c(p)|b|^p+p|a|^{p-2}a\cdot b
\]
where $c(p)> 0$. This yields

\[\rho^{\alpha}|\nabla_{\gamma}\rho|^t|\nabla\phi|^p\ge |\beta|^p
\rho^{\beta
p-p+\alpha}|\nabla_{\gamma}\rho|^{p+t}|\psi|^p+p|\beta|^{p-2}\beta
\rho^{\alpha +\beta
p+1-p}|\nabla_{\gamma}\rho|^{p+t-2}|\psi|^{p-2}\psi
\nabla\rho\cdot \nabla\psi.\] Integrating over the domain $\Omega$
gives

\begin{equation}\begin{aligned}\int_{\Omega}
\rho^{\alpha}|\nabla_{\gamma}\rho|^t|\nabla\phi|^pdx &\ge
|\beta|^p \int_{\Omega} \rho^{\beta
p-p+\alpha}|\nabla_{\gamma}\rho|^t|\psi|^pdz\\
&-p\int_{\Omega}|\beta|^{p-2}\beta \rho^{\alpha +\beta
p+1-p}|\nabla_{\gamma}\rho|^{p+t-2}|\psi|^{p-2}\psi
\nabla\rho\cdot \nabla\psi dz. \end{aligned}\end{equation}
Applying integration by parts to second integral on the right-hand
side of (3.10) and using the fact that
$\nabla_{\gamma}(|\nabla_{\gamma}\rho|)\cdot\nabla_{\gamma}\rho=0$
then we get \[\int_{\Omega}
\rho^{\alpha}|\nabla_{\gamma}\rho|^t|\nabla\phi|^pdx
\ge\Big(|\beta|^p-|\beta|^{p-2}\beta(\beta p-p+\alpha+Q)\Big)
\int_{\Omega} \rho^{\beta
p-p+\alpha}|\nabla_{\gamma}\rho|^{p+t}|\psi|^pdz. \] We now choose
$\beta=\frac{p-Q-\alpha}{p}$ to get the desired inequality
\begin{equation}
\int_{\Omega}
\rho^{\alpha}|\nabla_{\gamma}\rho|^t|\nabla\phi|^pdz\ge
\Big(\frac{Q+\alpha-p}{p}\Big)^p\int_{\Omega}\rho^{\alpha}|\nabla_{\gamma}\rho|^{t}\frac{|\nabla_{\gamma}\rho|^p}{\rho^
p}|\phi|^pdz.
\end{equation}
 Theorem (3.3) also holds for $1<p<2$
and in this case we use the following inequality
\[|a+b|^p-|a|^p\ge
c(p)\frac{|b|^2}{(|a|+|b|)^{2-p}}+p|a|^{p-2}a\cdot b
\]
where $c(p)> 0$ (see \cite {Lindqvist}). \qed
\end{theorem}
\medskip

We now have the following improved Hardy inequality which is
inspired by recent result of Abdellaoui, Colorado and Peral
\cite{ACP}. It is clear that if $\gamma=t=0$  then our result
recovers the inequality (1.4).

\begin{theorem}Let $\Omega\subset \mathbb{R}^{n}$ be a bounded domain with smooth
boundary which contains origin, $1<q<2$, $Q+\alpha-2>0$,
$Q=m+(1+\gamma)k$ and $\phi\in C_0^{\infty}(\Omega)$ then there
exists a positive constant $C=C(Q, q, \Omega)$ such that the
following inequality is valid
\begin{equation}\int_{\Omega}\rho^{\alpha}|\nabla_{\gamma}\rho|^{t}|\nabla_{\gamma}\phi|^2dz\ge
C_H\int_{\Omega}\rho^{\alpha}\frac{|\nabla_{\gamma}\rho|^{t+2}}{\rho^2}\phi^2dz+C\Big(\int_{\Omega}|\nabla_{\gamma}\phi|^q
\big(|\nabla_{\gamma}\rho|^t\rho^{\alpha}\big)^{\frac{q}{2}}dz\Big)^{2/q}\end{equation}
where $C_H=\big(\frac{Q+\alpha-2}{2}\big)^2$. \proof Let $\phi\in
C_0^{\infty}(\Omega) $ and $\psi=\rho^{\beta}$ where $\beta\in
\mathbb{R}\setminus \{0\}$. Then straightforward computation shows
that
\[|\nabla_{\gamma}\phi|^2-\nabla_{\gamma}(\frac{\phi^2}{\psi})\cdot\nabla_{\gamma}\psi=
\Big|\nabla_{\gamma}\phi-\frac{\phi}{\psi}\nabla_{\gamma}\psi\Big|^2.\]
Therefore
\[\begin{aligned}\int_{\Omega}\Big(|\nabla_{\gamma}\phi|^2-\nabla_{\gamma}(\frac{\phi^2}{\psi})\cdot\nabla_{\gamma}\psi\Big)\rho^{\alpha}|\nabla_{\gamma}\rho|^{t}dz &=
\int_{\Omega}\Big|\nabla_{\gamma}\phi-\frac{\phi}{\psi}\nabla_{\gamma}{\psi}\Big|^2\rho^{\alpha}|\nabla_{\gamma}\rho|^{t}dz\\
&\ge
c\Big(\int_{\Omega}\Big|\nabla_{\gamma}\phi-\frac{\phi}{\psi}\nabla_{\gamma}\psi\Big|^q
\rho^{\frac{q\alpha}{2}}|\nabla_{\gamma}\rho|^{\frac{qt}{2}}dz\Big)^{2/q}\end{aligned}\]
where we used the Jensen's inequality in the last step. Applying
integration by parts, we obtain

\[\begin{aligned}\int_{\Omega}\Big(|\nabla_{\gamma}\phi|^2-\nabla_{\gamma}(\frac{\phi^2}{\psi})\cdot\nabla_{\gamma}\psi\Big)\rho^{\alpha}|\nabla_{\gamma}\rho|^{t}dz&=
\int_{\Omega}|\nabla_{\gamma}\phi|^2\rho^{\alpha}|\nabla_{\gamma}\rho|^{t}dz\\ & +\frac{\beta}{\alpha+\beta}\int_{\Omega}\Big(\frac{\Delta_{\gamma}(\rho^{\alpha+\beta})}{\rho^{\beta}}\Big)|\nabla_{\gamma}\rho|^{t}\phi^2dz\\
&=\int_{\Omega}\rho^{\alpha}|\nabla_{\gamma}\rho|^{t}|\nabla_{\gamma}\phi|^2dz
\\ &+\beta(\alpha+\beta+Q-2)
\int_{\Omega}\rho^{\alpha}\frac{|\nabla_{\gamma}\rho|^{t+2}}{\rho^2}\phi^2dz.\end{aligned}\]
Therefore we have
\begin{equation}\begin{aligned}\int_{\Omega}\rho^{\alpha}|\nabla_{\gamma}\rho|^{t}|\nabla_{\gamma}\phi|^2dz&\ge
-\beta(\alpha+\beta+Q-2)\int_{\Omega}\rho^{\alpha}\frac{|\nabla_{\gamma}\rho|^{t+}2}{\rho^2}\phi^2dz\\
&+
c\Big(\int_{\Omega}\Big|\nabla_{\gamma}\phi-\frac{\phi}{\psi}\nabla_{\gamma}\psi\Big|^q
\rho^{\frac{q\alpha}{2}}|\nabla_{\gamma}\rho|^{\frac{qt}{2}}dz\Big)^{2/q}.\end{aligned}\end{equation}
We can use the following inequality which is valid for any  $w_1,
w_2\in \mathbb{R}^n $  and $1<q<2$
\begin{equation}c(q)|w_2|^q\ge |w_1+w_2|^q-|w_1|^q-q|w_1|^{q-2}\langle w_1,
w_2\rangle.\end{equation} Using the inequality (3.14), Young's
inequality and the weighted $L^p$-Hardy inequality (3.9), we get

\begin{equation}\int_{\Omega}\Big|\nabla_{\gamma}\phi-\frac{\phi}{\psi}\nabla_{\gamma}\psi\Big|^q
\rho^{\frac{q\alpha}{2}}|\nabla_{\gamma}\rho|^{\frac{qt}{2}}dz\ge
C \int _{\Omega}|\nabla_{\gamma}\phi|^q
\rho^{\frac{q\alpha}{2}}|\nabla_{\gamma}\rho|^{\frac{qt}{2}}dz
\end{equation}
where $C>0$.  Substituting (3.15)  into (3.13) then we obtain
\[\int_{\Omega}\rho^{\alpha}|\nabla_{\gamma}\rho|^t|\nabla_{\gamma}\phi|^2dz\ge -\beta(\alpha+\beta+Q-2)\int_{\Omega}\rho^{\alpha}\frac{|\nabla_{\gamma}\rho|^{t+2}}{\rho^2}\phi^2dz+
C\Big(\int_{\Omega}|\nabla_{\gamma}\phi|^q
\rho^{\frac{q\alpha}{2}}|\nabla_{\gamma}\rho|^{\frac{qt}{2}}dz\Big)^{2/q}.\]
Now choosing $\beta=\frac{2-\alpha-Q}{2}$ then we have the
following inequality
\[\int_{\Omega}\rho^{\alpha}|\nabla_{\gamma}\rho|^t|\nabla_{\gamma}\phi|^2dz\ge \Big(\frac{Q+\alpha-2}{2}\Big)^2\int_{\Omega}\rho^{\alpha}\frac{|\nabla_{\gamma}\rho|^{t+2}}{\rho^2}\phi^2dz+
C\Big(\int_{\Omega}|\nabla_{\gamma}\phi|^q
\rho^{\frac{q\alpha}{2}}|\nabla_{\gamma}\rho|^{\frac{qt}{2}}dz\Big)^{2/q}.\]
\end{theorem}
\endproof
\medskip

\section{Sharp Weighted Rellich-type inequalities}
 The main goal of this section is to find
sharp analogues of (1.6) and (1.7) for Baouendi-Grushin vector
fields. We then obtain their improved versions for bounded
domains. The proofs are mainly based on Hardy type inequalities.
The following is the first result of this section.

\begin{theorem}(Rellich type inequality I)
Let $\phi\in C_0^{\infty}(\mathbb{R}^{m+k}\setminus \{(0,0)\})$,
$Q=m+(1+\gamma)k$  and $\alpha>2$. Then the following inequality
is valid
\begin{equation}\begin{aligned}
\int_{\mathbb{R}^n}\frac{\rho^{\alpha}}{|\nabla_{\gamma }\rho|^2}
|\Delta_{\gamma} \phi|^2dz &\ge
\frac{(Q+\alpha-4)^2(Q-\alpha)^2}{16}
\int_{\mathbb{R}^n}\rho^{\alpha}\frac{|\nabla_{\gamma}
\rho|^2}{\rho^4}\phi^2dz.\end{aligned}\end{equation} Moreover, the
constant $\frac{(Q+\alpha-4)^2(Q-\alpha)^2}{16}$ is sharp.
\end{theorem}
\proof A straightforward computation shows that
\begin{equation}\Delta_{\gamma}\rho^{\alpha-2}=(Q+\alpha-4)(\alpha-2)\rho^{\alpha-4}|\nabla_{\gamma}\rho|^2.\end{equation} Multiplying both sides of (4.2) by $\phi^2$ and
integrating over $\mathbb{R}^n$,  we obtain

\[\int_{\mathbb{R}^n}\phi^2\Delta _{\gamma}\rho^{\alpha-2}dz=\int_{\mathbb{R}^n}\rho^{\alpha-2}(2\phi\Delta_{\gamma}\phi+2|\nabla_{\gamma}\phi|^2)dz.\]
Since
\[\int_{\mathbb{R}^n}\phi^2\Delta
_{\gamma}\rho^{\alpha-2}dz=(Q+\alpha-4)(\alpha-2)\int_{\mathbb{R}^n}\rho^{\alpha-4}|\nabla_{\gamma}\rho|^2\phi^2dz.\]
Therefore
\begin{equation}(Q+\alpha-4)(\alpha-2)\int_{\mathbb{R}^n}\rho^{\alpha-4}|\nabla_{\gamma}\rho|^2\phi^2dz-2\int_{\mathbb{R}^n}\rho^{\alpha-2}\phi\Delta_{\gamma}\phi
dx=
2\int_{\mathbb{R}^n}\rho^{\alpha-2}|\nabla_{\gamma}\phi|^2dz.\end{equation}
Applying the weighted Hardy inequality (3.9) to the right hand
side of (4.3), we get

\begin{equation}
\begin{aligned}-\int_{\mathbb{R}^n}\rho^{\alpha-2}\phi\Delta_{\gamma}\phi
dz&\ge
(\frac{Q+\alpha-4}{2})(\frac{Q-\alpha}{2})\int_{\mathbb{R}^n}\rho^{\alpha-4}|\nabla_{\gamma}\rho|^2\phi^2dz.
\end{aligned}\end{equation} We now apply the Cauchy-Schwarz
inequality to obtain

\begin{equation}
-\int_{\mathbb{R}^n}\rho^{\alpha-2}\phi\Delta_{\gamma}\phi dz\le
\Big(\int_{\mathbb{R}^n}\rho^{\alpha-4}|\nabla_{\gamma}\rho|^2\phi^2dz\Big)^{1/2}\Big(\int_{\mathbb{R}^n}\frac{\rho^{\alpha}}{|\nabla_{\gamma}\rho|
^2}|\Delta_{\gamma}\phi|^2dz\Big)^{1/2}.
\end{equation}
Substituting (4.5) into (4.4) yields the desired inequality

\begin{equation}\begin{aligned}
\int_{\mathbb{R}^n}\frac{\rho^{\alpha}}{|\nabla_{\gamma }\rho|^2}
|\Delta_{\gamma} \phi|^2dz &\ge
\frac{(Q+\alpha-4)^2(Q-\alpha)^2}{16}
\int_{\mathbb{R}^n}\rho^{\alpha}\frac{|\nabla_{\gamma}
\rho|^2}{\rho^4}\phi^2dz.\end{aligned}\end{equation}

It only remains to show that the constant $
C(Q,\alpha)=\frac{(Q+\alpha-4)^2(Q-\alpha)^2}{16}$ is the best
constant for the Rellich inequality (4.1), that is
\[\frac{(Q+\alpha-4)^2(Q-\alpha)^2}{16}=\inf\Big\{\frac{\int_{\mathbb{R}^n}\rho^{\alpha}\frac{|\Delta_{\gamma}f|^2}{|\nabla_{\gamma} \rho|^2}dz}{
\int_{\mathbb{R}^n}\rho^{\alpha}\frac{|\nabla_{\gamma}\rho|^2}{\rho^4}f^2dz},
f\in C_0^{\infty}(\mathbb{R}^n), f\neq 0\Big\}.\]  Given
$\epsilon>0$, take the radial function
\begin{equation}\phi_{\epsilon}(\rho)=\begin{cases}
 (\frac{Q+\alpha-4}{2}+\epsilon)\big(\rho-1\big)+1 &\quad\text{if}  \quad \rho\in [0,1],\\
\rho^{-(\frac{Q+\alpha-4}{2}+\epsilon)} &\quad \text{if} \quad
\rho>1,
\end{cases}
\end{equation}
where $\epsilon>0$. In the sequel we indicate
$B_1=\{\rho(z):\rho(z)\le 1\}$ $\rho$-ball centered at the origin
in $ \mathbb{R}^n$ with radius $1$.

By direct  computation we get
\begin{equation}\begin{aligned}\int_{B_{\rho}}\rho^{\alpha}\frac{|\Delta_{\gamma}\phi_{\epsilon}|^2}{|\nabla_{\gamma}
\rho|^2}dz&=\int_{B_1}\rho^{\alpha}\frac{|\Delta_{\gamma}\phi_{\epsilon}
|^2}{|\nabla_{\gamma} \rho|^2}dz+\int_{B_{\rho}\setminus
B_1}\rho^{\alpha}\frac{|\Delta_{\gamma}
\phi_{\epsilon}|^2}{|\nabla_{\gamma} \rho|^2}dz,\\
&=A(Q,\alpha, \epsilon)+B(Q, \alpha,
\epsilon)\int_{B_{\rho}\setminus
B_1}\rho^{-Q-2\epsilon}|\nabla_{\gamma}\rho|^2dz
\end{aligned}\end{equation}
where
\[ B(Q, \alpha,
\epsilon)=(\frac{Q+\alpha-4}{2}+\epsilon)^2(\frac{Q-\alpha}{2}-\epsilon)^2.\]
\begin{equation}
\begin{aligned}\int_{B_{\rho}}\rho^{\alpha}\frac{|\nabla_{\gamma}\rho|^2}{\rho^4}\phi^2dz&=\int_{B_1}
\rho^{\alpha}\frac{|\nabla_{\gamma}\rho|^2}{\rho^4}\phi^2dz+\int_{B_{\rho}\setminus
B_1}\rho^{\alpha}\frac{|\nabla_{\gamma}\rho|^2}{\rho^4}\phi^2dz  \\
&=C(Q,\alpha, \epsilon)+\int_{B_{\rho}\setminus B_1}
\rho^{-Q-2\epsilon}dz.
\end{aligned}
\end{equation}
Since $Q+\alpha-4>0$ then  $A(Q,\alpha,\epsilon)$,
 and $C(Q,\alpha,\epsilon)$ are bounded  and
we conclude by letting $\epsilon\longrightarrow 0$.\qed
\medskip

Using the same argument as above  and improved Hardy inequality
(3.1), we obtain the following improved Rellich type inequality.

\begin{theorem}
Let $\phi\in C_0^{\infty}(B_{\rho})$, $Q=m+(1+\gamma)k$  and
$4-Q<\alpha<Q$. Then the following inequality is valid
\begin{equation}\begin{aligned}
\int_{B_{\rho}}\frac{\rho^{\alpha}}{|\nabla_{\gamma }\rho|^2}
|\Delta_{\gamma} \phi|^2dz &\ge
\frac{(Q+\alpha-4)^2(Q-\alpha)^2}{16}
\int_{B_{\rho}}\rho^{\alpha}\frac{|\nabla_{\gamma}
\rho|^2}{\rho^4}\phi^2dz\\
&+\frac{(Q+\alpha-4)(Q-\alpha)}{2C^2r^2}\int_{B_{\rho}}\rho^{\alpha-2}\phi^2dz.\end{aligned}\end{equation}
\end{theorem}
\proof We have the following fact from (4.3):
\begin{equation}(Q+\alpha-4)(\alpha-2)\int_{B_{\rho}}\rho^{\alpha-4}|\nabla_{\gamma}\rho|^2\phi^2dz-2\int_{B_{\rho}}\rho^{\alpha-2}\phi\Delta_{\gamma}\phi
dx=
2\int_{B_{\rho}}\rho^{\alpha-2}|\nabla_{\gamma}\phi|^2dz.\end{equation}
Applying the improved Hardy inequality (3.1) on the right hand
side of (4.11), we get

\[\begin{aligned}&(Q+\alpha-4)(\alpha-2)\int_{B_{\rho}}\rho^{\alpha-4}|\nabla_{\gamma}\rho|^2\phi^2dz-2\int_{B_{\rho}}\rho^{\alpha-2}\phi\Delta_{\gamma}\phi
dz\\ &\ge 2(\frac{Q+\alpha-4}{2})^2
\int_{B_{\rho}}\rho^{\alpha-4}|\nabla_{\gamma}
\rho|^2\phi^2dz+\frac{2}{C^2r^2}\int_{B_{\rho}}\rho^{\alpha-2}\phi^2dz\end{aligned}\]
Now it is clear that,

\begin{equation}
\begin{aligned}-\int_{B_{\rho}}\rho^{\alpha-2}\phi\Delta_{\gamma}\phi
dz&\ge
(\frac{Q+\alpha-4}{2})(\frac{Q-\alpha}{2})\int_{B_{\rho}}\rho^{\alpha-4}|\nabla_{\gamma}\rho|^2\phi^2dz\\&+\frac{1}{C^2r^2}\int_{B_{\rho}}\rho^{\alpha-2}\phi^2dz.
\end{aligned}\end{equation} Next, we apply the Young's inequality to
the expression $-\int_{B_{\rho}}\rho^{\alpha-2}\phi\Delta\phi dz$
and we obtain
\begin{equation}-\int _{B_{\rho}}\rho^{\alpha-2}\phi\Delta_{\gamma}\phi dz\le \epsilon\int
_{B_{\rho}}\rho^{\alpha-4}|\nabla_{\gamma}\rho|^2\phi^2dz+\frac{1}{4\epsilon}\int
_{B_{\rho}}\rho^{\alpha}
\frac{|\Delta_{\gamma}\phi|^2}{|\nabla_{\gamma}\rho|^2}dz\end{equation}
where $\epsilon>0$.  Combining (4.13) and (4.12), we obtain

\[\int _{B_{\rho}}\rho^{\alpha}
\frac{|\Delta_{\gamma}\phi|^2}{|\nabla_{\gamma}\rho|^2}dz\ge
\big(-4\epsilon^2-(Q+\alpha-4)(Q-\alpha)\epsilon\big)\int
_{B_{\rho}}\rho^{\alpha-4}|\nabla_{\gamma}\rho|^2\phi^2dz+\frac{4\epsilon}{C^2r^2}\int_{B_{\rho}}\rho^{\alpha-2}\phi^2dz.\]
Note that the quadratic function
$-4\epsilon^2-(Q+\alpha-4)(Q-\alpha)\epsilon$ attains the maximum
for $\epsilon=\frac{(Q+\alpha-4)(Q-\alpha)}{8}$ and this maximum
is equal to $\frac{(Q+\alpha-4)^2(Q-\alpha)^2}{16}$. Therefore  we
obtain the desired inequality
\begin{equation}\begin{aligned}
\int_{B_{\rho}}\frac{\rho^{\alpha}}{|\nabla_{\gamma }\rho|^2}
|\Delta_{\gamma} \phi|^2dz &\ge
\frac{(Q+\alpha-4)^2(Q-\alpha)^2}{16}
\int_{B_{\rho}}\rho^{\alpha}\frac{|\nabla_{\gamma}
\rho|^2}{\rho^4}\phi^2dz\\
&+\frac{(Q+\alpha-4)(Q-\alpha)}{2C^2r^2}\int_{B_{\rho}}\rho^{\alpha-2}\phi^2dz.\end{aligned}\end{equation}\qed
\medskip

Arguing as above, and using the improved Hardy inequalities (3.2)
and  (3.4) we obtain the following Rellich type inequalities.

\begin{theorem}
Let $\phi\in C_0^{\infty}(\mathbb{R}^{m+k}\setminus \{(0,0)\})$,
$Q=m+(1+\gamma)k$  and $4-Q<\alpha<Q$. Then the following
inequality is valid
\begin{equation}\begin{aligned}
\int_{B_{\rho}}\frac{\rho^{\alpha}}{|\nabla_{\gamma }\rho|^2}
|\Delta_{\gamma} \phi|^2dz &\ge
\frac{(Q+\alpha-4)^2(Q-\alpha)^2}{16}
\int_{B_{\rho}}\rho^{\alpha-4}|\nabla_{\gamma}
\rho|^2\phi^2dz\\
&+\frac{(Q+\alpha-4)(Q-\alpha)}{8}\int_{B_{\rho}}\rho^{\alpha-4}|\nabla_{\gamma}\rho|^2\frac{\phi^2}{\ln
(\frac{r}{\rho})^2}dz.\end{aligned}\end{equation}
\end{theorem}
\begin{theorem}
Let $\phi\in C_0^{\infty}(\mathbb{R}^{m+k}\setminus \{(0,0)\})$,
$Q=m+(1+\gamma)k$ and $4-Q<\alpha<Q$. Then the following
inequality is valid
\begin{equation}\begin{aligned}
\int_{B_{\rho}}\frac{\rho^{\alpha}}{|\nabla_{\gamma }\rho|^2}
|\Delta_{\gamma} \phi|^2dz &\ge
\frac{(Q+\alpha-4)^2(Q-\alpha)^2}{16}
\int_{B_{\rho}}\rho^{\alpha}\frac{|\nabla_{\gamma}
\rho|^2}{\rho^4}\phi^2dz\\&+\frac{C(Q-\alpha)(Q+3\alpha-8)}{4}\big(\int_{\Omega}|\nabla_{\gamma}\phi|^q\rho^{\frac{q\alpha}{2}}dz\big)^{2/q}
\end{aligned}\end{equation} where $\Omega\subset \mathbb{R}^n$ is a bounded domain with smooth
boundary.
\end{theorem}
\bigskip

We now have the following Rellich type inequality that connects
first to second order derivatives. It is clear that if
$\alpha=\gamma=0$ then our result covers the inequality (1.7).
\begin{theorem}(Rellich type inequality II)
Let $\phi\in C_0^{\infty}(\mathbb{R}^{m+k}\setminus \{(0,0)\})$,
$Q=m+(1+\gamma)k$ and $2<\alpha<Q$. Then the following inequality
is valid
\begin{equation}\begin{aligned}
\int_{\mathbb{R}^n}\rho^{\alpha}\frac{|\Delta_{\gamma}
\phi|^2}{|\nabla_{\gamma }\rho|^2}dz \ge \frac{(Q-\alpha)^2}{4}
\int_{\mathbb{R}^n}\rho^{\alpha}\frac{|\nabla_{\gamma}
\phi|^2}{\rho^{2}}dz.
\end{aligned}\end{equation}
Furthermore, the constant
$C(Q,\alpha)=\big(\frac{Q-\alpha}{2}\big)^2$ is sharp.
\end{theorem}

\proof The proof of this theorem is similar to the proof Theorem
(4.1). Using the same argument as above, we have the following
from (4.3)
\begin{equation}-\int_{\mathbb{R}^n}\rho^{\alpha-2}\phi\Delta_{\gamma}\phi dx=\int_{\mathbb{R}^n}\rho^{\alpha-2}|\nabla_{\gamma}\phi|^2dz -\frac{(Q+\alpha-4)(\alpha-2)}{2}\int_{\mathbb{R}^n}\rho^{\alpha-4}|\nabla_{\gamma}\rho|^2\phi^2dz
.\end{equation} It is clear that $(Q+\alpha-4)(\alpha-2)>0$ and
using the  Hardy inequality (3.9) ($p=2, t=0$) we get

\begin{equation}-\int_{\mathbb{R}^n}\rho^{\alpha-2}\phi\Delta_{\gamma}\phi dz \ge
\frac{Q-\alpha}{Q+\alpha-4}\int_{\mathbb{R}^n}\rho^{\alpha-2}|\nabla_{\gamma}\phi|^2dz.
\end{equation} Let us apply Young's inequality to expression $-\int_{\mathbb{R}^n}\rho^{\alpha-2}\phi\Delta_{\gamma}\phi\,dz$ and we
obtain
\begin{equation}\begin{aligned}-\int_{\mathbb{R}^n}\rho^{\alpha-2}\phi\Delta_{\gamma}\phi
dz &\le \epsilon
\int_{\mathbb{R}^n}\rho^{\alpha-4}|\nabla_{\gamma}\rho|^2\phi^2dz+\frac{1}{4\epsilon}\int_{\mathbb{R}^n}\rho^{\alpha}\frac{|\Delta_{\gamma}\phi|^2}{|\nabla_{\gamma}\rho|^2}dz\\
& \le \epsilon\big(\frac{2}{Q+\alpha-4}\big)^2
\int_{\mathbb{R}^n}\rho^{\alpha-2}|\nabla_{\gamma}\phi|^2dz+\frac{1}{4\epsilon}\int_{\mathbb{R}^n}\rho^{\alpha}\frac{|\Delta_{\gamma}\phi|^2}{|\nabla_{\gamma}\rho|^2}dz\end{aligned}\end{equation}
where $\epsilon>0$ and will be chosen later. Substituting (4.20)
into (4.19) and rearranging terms, we get
\begin{equation}\begin{aligned}
\int_{\mathbb{R}^n}\rho^{\alpha}\frac{|\Delta_{\gamma}
\phi|^2}{|\nabla_{\gamma }\rho|^2}dz &\ge
\frac{-16\epsilon^2}{(Q+\alpha-4)^2}+4\big(\frac{Q-\alpha}{Q+\alpha-4}\big)\epsilon
\int_{\mathbb{R}^n}\rho^{\alpha}\frac{|\nabla_{\gamma}\phi|^2}{\rho^2}dz.\end{aligned}
\end{equation}
Choosing $\epsilon=\frac{1}{8}(Q-\alpha)(Q+\alpha-4)$ which yields
the  desired inequality
\begin{equation}
\int_{\mathbb{R}^n}\rho^{\alpha}\frac{|\Delta_{\gamma}
\phi|^2}{|\nabla_{\gamma }\rho|^2}dz \ge \frac{(Q-\alpha)^2}{4}
\int_{\mathbb{R}^n}\rho^{\alpha}\frac{|\nabla_{\gamma}
\phi|^2}{\rho^{2}}dz.
\end{equation}

To show that constant $ \big(\frac{Q-\alpha}{2}\big)^2$  is sharp,
we use the same sequence of functions (4.7) and we get
\[\frac{\int_{\mathbb{R}^n}\rho^{\alpha}\frac{|\Delta_{\gamma}\phi_{\epsilon}|^2}{|\nabla_{\gamma} \rho|^2}dz}{
\int_{\mathbb{R}^n}\rho^{\alpha}\frac{|\nabla_{\gamma}\phi_{\epsilon}|^2}{\rho^2}dz}\longrightarrow
\big(\frac{Q-\alpha}{2}\big)^2\] as $\epsilon\longrightarrow 0$.
\medskip

Now, using the same argument as above and improved Hardy
inequalities (3.1), (3.6) and (3.7) we obtain the following
improved Rellich type inequalities.
\begin{theorem}
Let $\phi\in C_0^{\infty}(B_{\rho})$, $Q=m+(1+\gamma)k$ and
$2<\alpha<Q$. Then the following inequality is valid
\begin{equation}\begin{aligned}
\int_{B_{\rho}}\rho^{\alpha}\frac{|\Delta_{\gamma}
\phi|^2}{|\nabla_{\gamma }\rho|^2}dz \ge \frac{(Q-\alpha)^2}{4}
\int_{B_{\rho}}\rho^{\alpha}\frac{|\nabla_{\gamma}
\phi|^2}{\rho^{2}}dz
+\frac{(Q-\alpha)(Q+3\alpha-8)}{4C^2r^2}\int_{B_{\rho}}\rho^{\alpha}\frac{\phi^2}{\rho^2}dz\end{aligned}\end{equation}
where $C>0$ and $r$ is the radius of the ball $B_{\rho}$.
\end{theorem}

\begin{theorem} Let $\Omega$ be a bounded domain with smooth
boundary $\partial \Omega$. Let $\phi\in C_0^{\infty}(\Omega)$,
$Q=m+(1+\gamma)k$ and $2<\alpha<Q$. Then the following inequality
is valid
\begin{equation}\begin{aligned}
\int_{\Omega}\rho^{\alpha}\frac{|\Delta_{\gamma}
\phi|^2}{|\nabla_{\gamma }\rho|^2}dz &\ge (\frac{Q-\alpha}{2})^2
\int_{\Omega}\rho^{\alpha}\frac{|\nabla_{\gamma}
\phi|^2}{\rho^{2}}dz+\tilde{C}\Big(\int_{\Omega}|\nabla_{\gamma}\phi|^q\rho^{\frac{q(\alpha-2)}{2}}dz\Big)^{2/q}\end{aligned}\end{equation}
where $\tilde{C}=\frac{C(Q-\alpha)(Q+3\alpha-8)}{4}$ and $C>0$.
\end{theorem}

\begin{theorem}
Let $\phi\in C_0^{\infty}(B_{\rho})$, $Q=m+(1+\gamma)k$ and
$2<\alpha<Q$. Then the following inequality is valid
\begin{equation}\begin{aligned}
\int_{B_{\rho}}\rho^{\alpha}\frac{|\Delta_{\gamma}
\phi|^2}{|\nabla_{\gamma }\rho|^2}dz \ge \frac{(Q-\alpha)^2}{4}
\int_{B_{\rho}}\rho^{\alpha}\frac{|\nabla_{\gamma}
\phi|^2}{\rho^{2}}dz
+C(Q,\alpha)\int_{B_{\rho}}\rho^{\alpha-4}|\nabla_{\gamma}\rho|^2\frac{\phi^2}{(\ln
\frac{r}{\rho})^2}dz\end{aligned}\end{equation} where
$C(Q,\alpha)=\frac{(Q-\alpha)(Q+3\alpha-8)}{16}$.
\end{theorem}
\bibliographystyle{amsalpha}

\begin{thebibliography}{A}

\bibitem [ACP] {ACP}  B. Abdellaoui, D. Colorado, I. Peral, \textit{ Some improved Caffarelli-Kohn-Nirenberg inequalities},
 Calc. Var. Partial Differential Equations {\bf 23}  (2005), no. 3,
 327-345.
\bibitem [AR] {Adimurthi} Adimurthi N. Chaudhuri and M. Ramaswamy, \textit{An improved Hardy–Sobolev inequality and its applications},
 Proc. Amer. Math. Soc. {\bf 130} (2002), pp. 489–505.

\bibitem[BG]{Baras-Goldstein} P. Baras and J. A. Goldstein, \textit
{The heat equation with a singular potential}, Trans. Amer. Math.
Soc. {\bf 284} (1984), 121-139.

\bibitem[B] {Baouendi} M. S. Baouendi, \textit{Sur une classe d'op\'erateurs elliptiques d\'eg\'en\'er\'es}, Bull. Soc. Math. France {\bf 95}, 45-87 (1967).

\bibitem[Be] {Beliche} Bella\" {\i}che, Andr\' {e}.
\emph{The Tangent Space in Sub-Riemannian Geometry}, 1--78, Progr.
Math., 144, Birkhäuser, Basel, 1996.


\bibitem [BV]{Brezis-Vazquez} H. Brezis and J. L. V\'azquez, \textit{Blow-up solutions of some nonlinear elliptic
problems}, Rev. Mat. Univ. Complutense Madrid \textbf{10} (1997),
443-469.

\bibitem [CM]{Cabre-Martel} X. Cabr\'e and Y. Martel, \textit{Existence versus explosion instantanée pour des équations de la chaleur linéaires avec potentiel singulier}
, C. R. Acad. Sci. Paris Sér. I. Math., \textbf{329} (1999),
973-978.

\bibitem [D] {D'Ambrosio} L. D'Ambrosio, \textit{Hardy inequalities related to Grushin type operators},
 Proc. Amer. Math. Soc. {\bf 132} (2004), no. 3, 725-734.

\bibitem [DH] {Davies-Hinz} E. B. Davies, and A. M. Hinz, \textit{Explicit constants for Rellich inequalities in $L\sb p(\Omega)$},
 Math. Z. \textbf {227} (1998), no. 3, 511-523.

\bibitem [FGW] {Franchi} B. Franchi, C. E. Gutiérrez, and R. L. Wheeden, Weighted
Sobolev-Poincaré inequalities for Grushin type operators, Comm.
Partial Differential Equations 19, 523-604 (1994).

\bibitem [FGaW] {FGaW} B. Franchi, S. Gallot and R. L. Wheeden \textit{Sobolev and isoperimetric inequalities for degenerate metrics},
 Math. Ann. {\bf 300} (1994), no. 4, 557-571.

 \bibitem [FT] {FT} S. Filippas and  A. Tertikas \textit{Optimizing Improved Hardy Inequalities},
 J. Funct. Anal. {\bf 192}, (2002),186-233.

\bibitem [G] {Garofalo} N. Garofalo,  \textit{Unique continuation for a class of elliptic
operators which degenerate on a manifold of arbitrary
codimension},  J. Differential Equations {\bf 104} (1993), no. 1,
117-146.
\bibitem [GP]{Garcia-Peral}
J. Garcia Azorero and I. Peral, \textit{Hardy inequalities and
some critical elliptic and parabolic problems}, J. Diff.
Equations, {\bf 144} (1998), 441-476.
\bibitem [GK]{Goldstein-Kombe} J. A. Goldstein and I. Kombe,
\textit {Nonlinear parabolic differential equations with the
singular lower order term}, Adv. Differential Equations {\bf 10}
(2003), 1153-1192.
\bibitem [G1] {Grushin1} V. Grushin, \textit{A certain class of hypoelliptic operators}, Math. USSR-Sb. 12, No. 3, 458-476 (1970)

\bibitem [G2] {Grushin2} V. Grushin, \textit{A certain class of elliptic pseudodifferential operators that are degenerate on a submanifold}, Mat. Sb. 84, 163-195 (1971).

\bibitem [K1]{Kombe1} I. Kombe, \textit{Nonlinear degenerate parabolic equations for Baouendi-Grushin operators},
Math. Nachr. {\bf 279}  (2006), no. 7, 756-773.

\bibitem [K2]{Kombe2} I. Kombe, \textit{ Hardy, Rellich and Uncertainty principle inequalities on Carnot
Groups}, preprint.

\bibitem [K3]{Kombe3} I. Kombe, \textit{ Sharp Hardy and Rellich type inequalites with remainders for the Greiner vector fields}, preprint.

\bibitem [K\"O]{Kombe-Ozaydin} I. Kombe and M. \"Ozaydin  \textit{ Improved Hardy and Rellich inequalities on Riemannian manifolds},
preprint.


\bibitem [L] {Lindqvist} P. Lindqvist, \textit {On the equation}
$\text{div}(|\nabla u|^{p-2}\nabla u)+\lambda|u|^{p-2}u=0$, Proc.
Amer. Math. Soc. \textbf(109) (1990), 157-164.

\bibitem [Lu] {Lu} G. Lu, \textit{ Weighted Poincaré and Sobolev inequalities for vector fields satisfying Hörmander's condition and applications},
 Rev. Mat. Iberoamericana {\bf 8} (1992), no. 3, 367-439.

\bibitem [M] {Monti} R. Monti \textit{Sobolev inequalities for weighted gradients},
 Comm. Partial Differential Equations {\bf 31} (2006), no. 10-12,
 1479-1504.

\bibitem [PV]{Peral-Vazquez} I. Peral and J. L. V\'azquez,
\textit{ On the stability or instability of the singular solution
of the semilinear heat equation with exponential reaction term},
Arch. Rational Mech. Anal. {\bf 129} (1995),  201-224.

\bibitem [TZ] {TZ} A. Tertikas and  N. Zographopoulos,  \textit{Best constants in the Hardy-Rellich Inequalities and Related
Improvements},  Adv. Math. {\bf 209},2 (2007), 407-459.

\bibitem [VZ]{Vazquez-Zuazua} J. L. V\'azquez and E. Zuazua, \textit{The Hardy constant and the asymptotic behaviour of the heat equation with an inverse-square potential}
, J. Funct. Anal. \textbf{173} (2000), 103-153.

\bibitem [WW] {Wang-Willem} Z.-Q. Wang and  M. Willem, \textit{Caffarelli-Kohn-Nirenberg
inequalities with remainder terms},  J. Funct. Anal. {\bf 203}
(2003), no. 2, 550-568.

\end{thebibliography}

\end{document}